\pdfoutput=1
\documentclass[11pt]{article}

\usepackage{amssymb}

\usepackage{alltt,times,pifont,amssymb,url,graphicx} 
\usepackage{latexsym} 

\def\SAH{\cite{solovay-et-al09}}
\def\SAHr#1{\cite[#1]{solovay-et-al09}}

\def\Hide#1{\relax}

\DeclareSymbolFont{AMSb}{U}{msb}{m}{n}
\DeclareSymbolFontAlphabet{\mathbb}{AMSb}

\def\Func#1{{\mathsf{#1}}}

\def\VS#1{\Func{VS}^{#1}}

\def\LL{{\cal L}}

\def\LN{\LL_{N}}
\def\LNA{\LL^{+}_{N}}

\def\sL{\mathbb{L}}
\def\sLO{\sL_0}
\def\sLI{\sL_1}

\def\FRAC#1#2{\frac{{\displaystyle #1}}{{\displaystyle #2}}}
\def\Frac#1#2{\frac{{\textstyle #1}}{{\textstyle #2}}}

\def\Sin#1{\Func{sin}(#1)}

\def\pDef{\Func{Def}}

\def\pG#1{\Func{G}(#1)}
\def\pGamma#1{{\sf\Gamma}\left(#1\right)}

\def\pMultX#1{\Func{Mult}(#1)}
\def\pN#1{\Func{N}(#1)}
\def\pNNMult#1{\Func{Mult_{{\ge}0}}(#1)}

\def\pOK#1{\Func{OK}(#1)}

\def\pPar#1{\Func{Par}(#1)}
\def\pPeriodic{\Func{Periodic}}
\def\pPi#1{{\sf\Pi}\left(#1\right)}

\def\pRotund#1{\Func{Rotund}(#1)}

\def\pSD#1{\Func{SD}(#1)}
\def\pSIN#1{\Func{SIN}(#1)}

\def\pW#1{\Func{W}(#1)}

\def\fAE{\Al\Ex}
\def\fEA{\Ex\Al}
\def\fAIA{{\Al}{\Imp}{\Al}}

\def\IsDef{\mathrel{{:}{=}}}
\def\Al{\forall}
\def\Ex{\exists}
\def\Imp{\Rightarrow}

\def\Iff{\Leftrightarrow}
\def\And{\land}
\def\Or{\lor}
\def\Not{\lnot}

\def\ST{\mathrel{|}}

\def\Va{\mathbf{a}}
\def\Vb{\mathbf{b}}

\def\Ve{\mathbf{e}}

\def\Vp{\mathbf{p}}
\def\Vq{\mathbf{q}}
\def\Vr{\mathbf{r}}

\def\Vu{\mathbf{u}}
\def\Vv{\mathbf{v}}
\def\Vw{\mathbf{w}}
\def\Vx{\mathbf{x}}
\def\Vy{\mathbf{y}}
\def\Vz{\mathbf{z}}

\def\VA{\mathbf{A}}

\def\VS{\mathbf{S}}
\def\VT{\mathbf{T}}
\def\VU{\mathbf{U}}
\def\VV{\mathbf{V}}

\def\VX{\mathbf{X}}

\def\VZ{\mathbf{Z}}

\def\VZero{\mathbf{0}}
\def\VOne{\mathbf{1}}
\def\VTwo{\mathbf{2}}

\def\VFour{\mathbf{4}}

\def\Norm#1{\left|\left|{#1}\right|\right|}

\def\Diff{\mathop{\backslash}}

\def\Bndry#1{\partial{#1}}
\def\Int#1{{#1}^{\circ}}

\newtheorem{Theorem}{Theorem}
\newtheorem{Lemma}[Theorem]{Lemma}
\newtheorem{Corollary}[Theorem]{Corollary}

\def\Proof{\par \noindent{\bf Proof: }}
\def\Done{\hfill\rule{0.5em}{0.5em}}

\newcommand{\nat}{\mbox{$\protect\mathbb N$}}

\newcommand{\real}{\mbox{$\protect\mathbb R$}}

\newcommand{\realgtz}{\real_{{>}0}}

\newcommand{\spot}{{\cdot}}

\newcommand{\all}[1]{\forall #1 \spot\:}
\newcommand{\ex}[1]{\exists #1 \spot\:}

\newcommand{\Ands}{\bigwedge}

\newcommand{\BEQ}{\mbox{\raise4pt\hbox{$\ulcorner$}}}

\newcommand{\EEQ}{\mbox{\raise4pt\hbox{$\urcorner$}}}

\newcommand{\BA}{\begin{array}[t]{l}}
\newcommand{\EA}{\end{array}}

\newcommand{\To}{\rightarrow}

\makeatletter
\def\imod#1{\allowbreak\mkern10mu({\operator@font mod}\,\,#1)}
\makeatother

\newcommand\HOLSpacing{12pt}

\newlength{\hsbw}
\title{Undecidability for the additive $\fAIA$ fragment of the
theory of normed spaces.}
\author{R. D. Arthan \\
{\small
Lemma 1. Ltd., 2nd Floor, 31A Chain St, Reading, RG1 2HX, UK}  \\ \&
{\small
Department of Computer Science, Queen Mary,
University of London, London, E1 4NS, UK.} \\
{\small\tt rda@lemma-one.com}}
\date{6th February 2010\thanks{With formatting corrections dated 9th April 2011}}

\begin{document}
\maketitle

\begin{abstract}
An $\fAIA$ formula is one of the form $A \Imp B$ where $A$ and $B$ are purely
universal.  Up to a simple reduction $\fAIA$ formula are both $\fEA$ and $\fAE$.
In an earlier paper Solovay, Harrison and I proved the undecidability of validity for the $\fAIA$
fragment of a two-sorted first-order language $\LN$ for normed vector
spaces.   In this note we find that validity remains undecidable
for $\fAIA$ sentences in the additive sublanguage $\LNA$, i.e., when multiplication is disallowed.
\end{abstract}

Sections 7 and 8 of {\SAH} consider the
decision problem for a two-sorted language $\LN$ for normed real
vector spaces.  Decision procedures are given for the universal and existential
fragments and validity for $\fAIA$ sentences is shown to be undecidable.
In this note I sharpen these results by proving the undecidability of validity
for $\fAIA$ sentences in the additive sublanguage $\LNA$.  The proof for $\LN$
used a 2-dimensional space whose unit circle encodes the graph of the sine
function.  The proof given for $\LNA$ here is based on a modification of this
space such that scalar-scalar multiplication becomes definable and
scalar-vector multiplication is no longer needed.  Expressing multiplication
geometrically introduces some notational complexity, but the overall structure
of the argument remains the same.

We will need an additive way to assert that two
vectors in a normed space $V$ point in the same direction. 
For $\VZero \not = \Vv, \Vw \in V$ define $\Va(\Vv, \Vw)$ by:
\[
\begin{array}{lcl}
\Va(\Vv, \Vw) &\IsDef&
         \frac{\Norm{\Vv}}{\Norm{\Vv}+\Norm{\Vw}}\Vv +
         \frac{\Norm{\Vw}}{\Norm{\Vv}+\Norm{\Vw}}
                \frac{\Norm{\Vv}}{\Norm{\Vw}} \Vw
\end{array}
\]
\noindent so that $\Va(\Vv, \Vw)$ is a proper convex
combination of $\Vv$ and $\frac{\Norm{\Vv}}{\Norm{\Vw}}\Vw$
and we have:
\[
\begin{array}{rcccccl}
\Norm{\Va(\Vv, \Vw)} &=&
   \Norm{\frac{\Norm{\Vv}}{\Norm{\Vv}+\Norm{\Vw}}(\Vv + \Vw)} &=&
   \frac{\Norm{\Vv}\cdot\Norm{\Vv+\Vw}}{\Norm{\Vv}+\Norm{\Vw}} &\le&
\Norm{\Vv}
\end{array}
\]
\noindent with equality iff $\Norm{\Vv+\Vw} = \Norm{\Vv} + \Norm{\Vw}$ in
which case the line segment $[\Vv, \frac{\Norm{\Vv}}{\Norm{\Vw}}\Vw]$
is contained in the $V$-circle $S_{\Norm{\Vv}} = \{\Vu \ST \Norm{\Vu} = \Norm{\Vv}\}$.
It follows that if we define:
\[
\begin{array}{lcl}
\pSD{\Vv, \Vw} &\IsDef& \Norm{\Vv + \Vw} = \Norm{\Vv} + \Norm{\Vw},
\end{array}
\]
\noindent then for $\Vv, \Vw \not= \VZero$ with
$\frac{\Norm{\Vv}}{\Norm{\Vw}}\Vw \not = \Vv$, $\pSD{\Vv, \Vw}$ holds
iff the line segment $[\Vv, \frac{\Norm{\Vv}}{\Norm{\Vw}}\Vw]$ is 
contained in $S_{\Norm{\Vv}}$.
A {\em rotund point} of a normed space is a point $\Vv$ that is
not an end-point of any proper line segment contained in $S_{\Norm{\Vv}}$:
\[
\begin{array}{lcl}
\pRotund{\Vv} &\IsDef& \all{\Vu} \Norm{\Vu} = \Norm{\Vv} =
\Norm{(\Vu+\Vv)/2} \Imp \Vu = \Vv
\end{array}
\]
\noindent
By the remarks above, if $\Vw$ is a non-zero rotund point of $V$, then
$\pSD{\Vv, \Vw}$ holds iff $\Vv = \frac{\Norm{\Vv}}{\Norm{\Vw}}\Vw$.
Thus we need a space with an adequate supply of rotund points.

We write $S_V(\Vp, r)$ for
the circle with centre $\Vp$ and radius $r$ in the normed space $V$.
We will use the following generalisation of the fact that
two distinct euclidean circles meet in at most two points.
\begin{Lemma}\label{lma:two-point-lemma}
Let $V$ be a two-dimensional normed space and let $R$ be
the set of rotund points of $V$.
Let $\Va, \Vb, \Vp, \Vq \in V$ and let $r, s \in \real$ be such
that the following conditions hold:
$\Va \not= \Vb$;
$\Vp \not= \Vq$;
$\{\Va, \Vb\} \subseteq S_V(\Vp, r) \cap S_V(\Vq, s)$;
$\{\Va - \Vp, \Va - \Vq\} \cap R \not= \emptyset$ and
$\{\Vb - \Vp, \Vb - \Vq\} \cap R \not= \emptyset$.
Then $S_V(\Vp, r) \cap S_V(\Vq, s) = \{\Va, \Vb\}$.
\end{Lemma}
\Proof See Appendix~\ref{sec:intersecting-circles}. \Done

The proof that we wish to adapt \SAHr{theorem~45}
involves a 2-dimensional space $\sLO$ in which the
north-west quadrant of the unit circle is the graph of a
function $\gamma : (-1, 0) \rightarrow (0, 1)$ defined
as follows:
\[
\begin{array}{rcl}
\gamma(x) &=& \Frac{g(\frac{x+1}{-x})}{1 + g(\frac{x+1}{-x})}.
\end{array}
\]
where $g(s) = 2s + s^2 + \frac{1}{M}\Sin{s}$ for a certain positive integer $M$.
In the proof of \SAHr{lemma~43}, it is shown
that $\gamma''$ is defined and negative throughout $(-1, 0)$, which implies
that the graph of $\gamma$ contains no proper line segments. This implies that every point in
the interior of the north-west quadrant of $\sLO$ is rotund.  However, the
points $-\Ve_1$ and $\Ve_2$ are not rotund in $\sLO$.  To rectify this, we must
redesign the north-east quadrant.
To describe the redesigned norm we will use the $\sLO$-norm, but
only for vectors in the interiors of the north-west and
south-east quadrants, where the norm is determined by $\gamma$.

\begin{figure}
\begin{center}
\includegraphics[angle=0,scale=0.4]{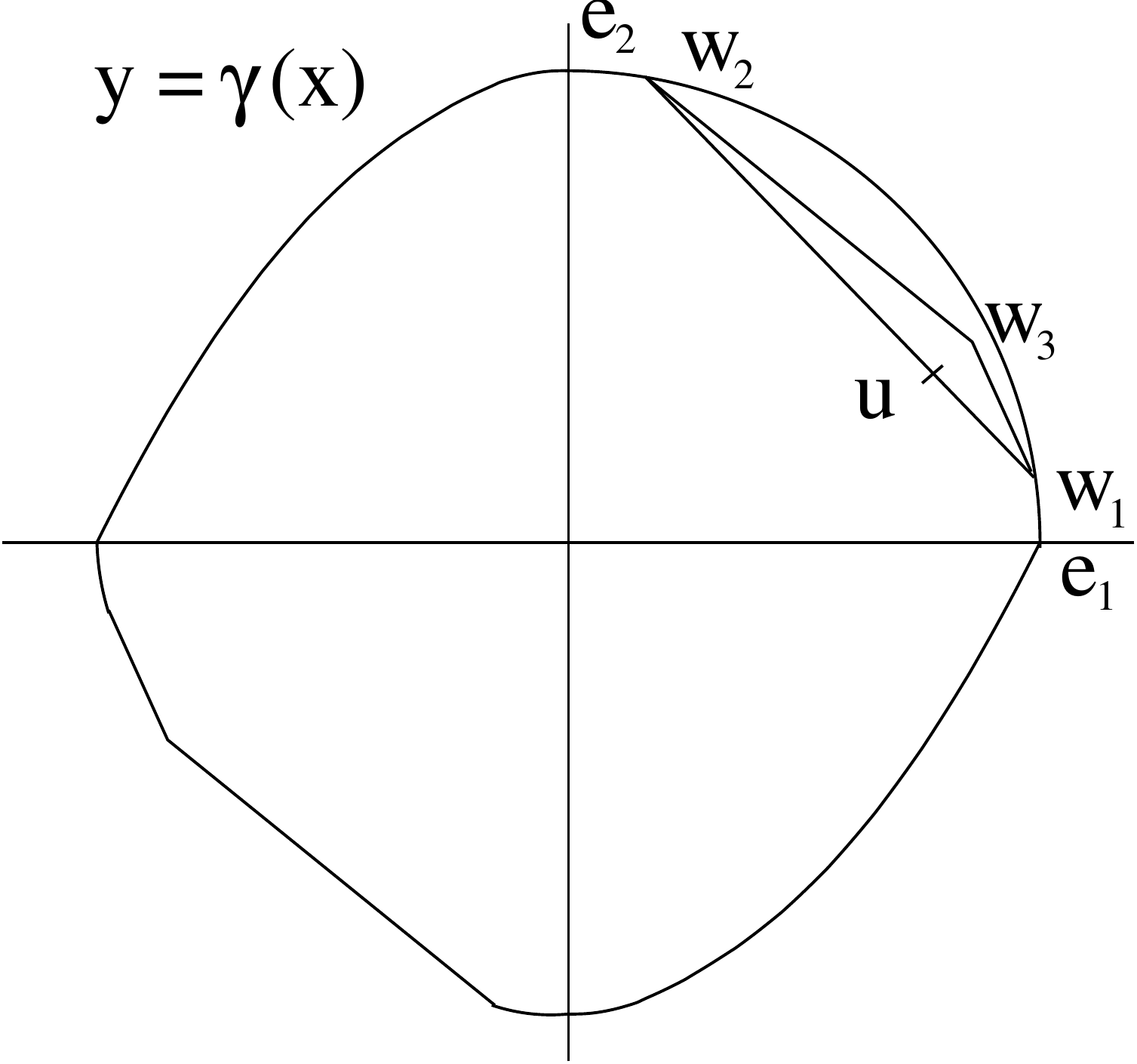}
\caption{The normed space $\sLI$}
\label{fig:l1}
\end{center}
\end{figure}

Let $\Vw_i$, $i = 1, 2,$ be the euclidean unit vectors in the north-east
quadrant such that $\Norm{\Ve_i - \Vw_i}_{\sLO} = q$ where the
rational constant $q$ is chosen so that $0 < q < \frac{1}{4}$ and
$d = \Norm{\Vw_1 - \Vw_2}_{\sLO} > \frac{3}{4}$ (this is possible
since $\Norm{\Ve_1 - \Ve_2}_{\sLO} > 1$.).
For $i =
1, 2$ and $r > 0$, let $S_i(r) = S_{\sLO}(\Vw_i, r)$.
$S_1(\frac{d}{3})$ and $S_2(\frac{2d}{3})$ meet at the point $\Vu =
\frac{2}{3}\Vw_1 + \frac{1}{3}\Vw_2$ and for small enough $r > \frac{d}{3}$, $S_1(r)
\cap S_2(2r)$ comprises two points, one on each side of the line from
$\Vw_1$ to $\Vw_2$, and these converge to $\Vu$ as $r$ tends to $\frac{d}{3}$
(to see this, show that for small enough $r > \frac{d}{3}$, one has {\em(a)}
that for any $\Vx \in S_1(r) \cap S_2(2r)$ the vectors $\Vx - \Vw_1$
and $\Vx - \Vw_2$ are both rotund in $\sLO$ and {\em(b)} that
$S_1(r) \cap S_2(2r)$ contains at least two points;
now apply lemma~\ref{lma:two-point-lemma} to conclude
that $S_1(r) \cap S_2(2r)$ contains exactly two points for small $r > \frac{d}{3}$
and note that if $\Vx$ is a point of accumulation of the set
$\bigcup_{n=1}^{\infty} (S_1(\frac{d}{3}+\frac{1}{n}) \cap S_2(\frac{2d}{3}+\frac{2}{n}))$,
then $\Vx \in S_1(\frac{d}{3}) \cap S_2(\frac{2d}{3}) = \{\Vu\}$.).

Writing $\Norm{\Vu}_e$ for the euclidean norm, we have $\Norm{\Vu}_e <
1$, and so we can choose a rational constant $r > \frac{d}{3} \ge \frac{1}{4}$
and a
point $\Vw_3 \in S_1(r) \cap S_2(2r)$ with $\Norm{\Vw_3}_{e} < 1$ lying
to the north-east of the line segment $[\Vw_1, \Vw_2]$.  Let $U$ be the
set comprising: the point $-\Ve_1$, the graph of $\gamma$, the arc of
the euclidean unit circle from $\Ve_2$ to $\Vw_2$, the line segment
$[\Vw_2, \Vw_3]$, the line segment $[\Vw_3, \Vw_1]$ and the arc of the
euclidean unit circle from $\Vw_3$ to $\Ve_1$ (see Figure~\ref{fig:l1}.).
Using the fact that $\gamma''$ is defined and negative throughout $(-1,
0)$, it is straightforward to verify that $U$ is the graph of a
concave function, so $S = U \cup -U$ is indeed the unit
circle of a norm, and we define $\sLI$ to be $\real^2$ equipped
with this norm. 
The maximal line segments in $S$  
 are the segments $\pm[\Vw_1, \Vw_3]$ and $\pm[\Vw_3, \Vw_2]$
and so all points of the north-west quadrant are
rotund in $\sLI$ including $-\Ve_1$ and $\Ve_2$.  These  maximal line
segments have rational $\sLI$-lengths $r$ and $2r$, since the $\sLO$ and $\sLI$ norms agree in the north-west quadrant.

We will give replacements in $\LNA$ for the various formulas
of $\LN$ used in the proof of {\SAHr{theorem~45}}.
We begin with the replacement for $\pW{\Vp, \Vq, \Vr}$ which now
has to characterize a configuration of 5 points rather than 3:
\[
\begin{array}{lcl}
\pW{\Vp_1, \Vp_2, \Vu_1, \Vu_2, \Vu_3} &\IsDef&
    \Norm{\Vp_1} = \Norm{\Vp_2} = \Norm{\Vu_1} = \Norm{\Vu_2} = \Norm{\Vu_3} = 1 \\
&\And&
    \Norm{(\Vu_1 + \Vu_3)/2} = \Norm{(\Vu_2 + \Vu_3)/2} = 1 \\
&\And&
    \Norm{(\Vu_1 + \Vu_2)/2} < 1 \\
&\And&
    \Norm{\Vu_1 - \Vu_3} = r \And \Norm{\Vu_3 - \Vu_2} = 2r \\
&\And&
    \Norm{\Vp_1 - \Vu_1} =  \Norm{\Vp_2 - \Vu_2} = q \\
&\And&
    \Norm{(\Vp_1 + \Vu_1)/2} < 1 \And \Norm{(\Vp_2 + \Vu_2)/2} < 1 \\
\end{array}
\]
\noindent
Here $q$ and $r$ are the rational constants chosen in the construction
of $\sLI$.
The conditions on $\Vu_1$, $\Vu_2$ and $\Vu_3$ imply that the edges $[\Vu_1,
\Vu_3]$ and $[\Vu_2, \Vu_3]$ of the triangle $\bigtriangleup\Vu_1\Vu_2\Vu_3$
are contained in the unit circle and have lengths $r$ and $2r$ respectively.
For $i = 1, 2$, the conditions on $\Vp_i$ imply that $\Vp_i \in S$ and
$\Norm{\Vp_i, \Vu_i} = q$ but $\Vp_i \not\in [\Vu_i, \Vu_3]$.  So $\pW{\Vp_1,
\Vp_2, \Vu_1, \Vu_2, \Vu_3}$ holds in $\sLI$ iff $(\Vp_1, \Vp_2, \Vu_1, \Vu_2,
\Vu_3) = \pm (\Ve_1, \Ve_2, \Vw_1, \Vw_2, \Vw_3)$ (see this using
lemma~\ref{lma:two-point-lemma} to show that $S_{\sLI}(\Vw_i, q) \cap S = \{\Ve_i, \Vx\}$
say, where $\Vx \in [\Vw_i, \Vw_3]$, because $q < \frac{1}{4} < r$, while
$\Ve_i \not\in [\Vw_i, \Vw_3]$.).

To give the replacement for the formula $\pDef$ it is convenient
to define the following formula $\pPar{\Vv, \Vw}$ abstracting the notion that $\Vv$ and $\Vw$ are parallel:
\[
\begin{array}{lcl}
\pPar{\Vv, \Vw} &\IsDef& \Vv \not= \VZero \And \Vw \not= \VZero \And 
    (\pSD{\Vv, \Vw} \Or \pSD{\Vv, -\Vw})
\end{array}
\]
For rotund $\Vv$, $\pPar{\Vv, \Vw}$ holds iff $\Vv$ and $\Vw$ span the
same line through the origin.
Even if  $\Vv, \Vw \not= \VZero$ are not rotund, if $\pPar{\Vv, \Vw}$ fails to
hold, then $\Vv$ and $\Vw$ must be linearly independent.  Now we give the
formula $\pDef$ using vector variables $\Vx$ and $\Vy$ constrained
to be parallel to the coordinate axes to stand in for the
vectors $x\Ve_1$ and $y\Ve_2$ of the original.
\[
\begin{array}{lcl}
\pDef &\IsDef&
    (\Norm{\Ve_1} = \Norm{\Ve_2} = 1 \And \Not\pPar{\Ve_1, \Ve_2}) \\
&\And&
    (\pPar{\Vx, \Ve_1} \And \pPar{\Vy, \Ve_2} \And
                       \Norm{\Vx + \Vy} = 1 \Imp 
    \Norm{\Vx} < 1 \And \Norm{\Vy} < 1) \\
&\And&
    (\pSD{\Vx, -\Ve_1}  \And \pSD{\Vy, \Ve_2} \And
         \Norm{\Vz} = \Norm{\Vx + \Vy} = \Norm{(\Vx + \Vy +\Vz)/2} \\
&& \quad\quad\quad\quad \Imp \Vz = \Vx + \Vy)
\end{array}
\]
As before, the first conjunct above implies that $\Ve_1$ and $\Ve_2$ are
linearly independent unit vectors and, provided $\Ve_1$ and $\Ve_2$ are rotund
points, the second conjunct implies that in
the plane spanned by $\Ve_1$ and $\Ve_2$, the unit disc is contained
in the square with diagonal $[-\Ve_1-\Ve_2, \Ve_1 + \Ve_2]$ and meets
the edges of that square in the points $\pm\Ve_1$ and $\pm\Ve_2$..
The third conjunct is new: for $\Vx = x\Ve_1$
and $\Vy = y\Ve_2$ with $x \le 0 \le y$, it holds for all $\Vz$
iff $\Vx + \Vy$ is rotund.

To give the replacement for the formula $\pG{s, t}$ we need to represent
scalar-scalar multiplication. It will be convenient to represent a real
number $s$ by the pair of vectors $(-s\Ve_1, s\Ve_2)$.
We use bold uppercase letters $\VS$, $\VT$ etc., as abbreviations for
pairs of vector variables and write $\VS.1$ and $\VS.2$ for the two
components of $\VS$.  Thus the following predicate defines the pairs
that represent real numbers:
\[
\begin{array}{lcl}
\pOK{\VS} &\IsDef& \Norm{\VS.1} = \Norm{\VS.2} \\
&\And&
    \left(\strut(\pSD{\VS.1, -\Ve_1} \And \pSD{\VS.2, \Ve_2}) \Or (\pSD{\VS.1, \Ve_1} \And \pSD{\VS.2, -\Ve_2})\right) \\
\end{array}
\]
\noindent
We introduce some syntactic abbreviations for working with
pairs of vectors, while taking care at each point of
use that the abbreviations can
be expanded out to give a genuine formula of $\LNA$.
We add, negate and scale pairs of vector expressions component-wise:
$\VS+\VT \IsDef (\VS.1 + \VT.1, \VS.2 + \VT.2)$, $-\VS \IsDef (-\VS.1,
-\VS.2)$ and $x\VS \IsDef (x\VS.1, x\VS.2)$ (where $x$ will always be a
rational constant so that we remain in $\LNA$).  We write $\VZero$,
$\VOne$, $\VTwo, \ldots$ for the pairs $(-i\Ve_1, i\Ve_2)$, $i = 0, 1,
2, \ldots$  We write $\VS = \VT$ for $\VS.2 =
\VT.2$, $\VS \ge \VT$ or $\VT \le \VS$ for $\pSD{(\VS - \VT).2, \Ve_2}$
and $\VS > \VT$ or $\VT < \VS$ for $\VS \ge \VT \And \VS \not= \VT$.

\begin{figure}
\begin{center}
\includegraphics[angle=0,scale=0.7]{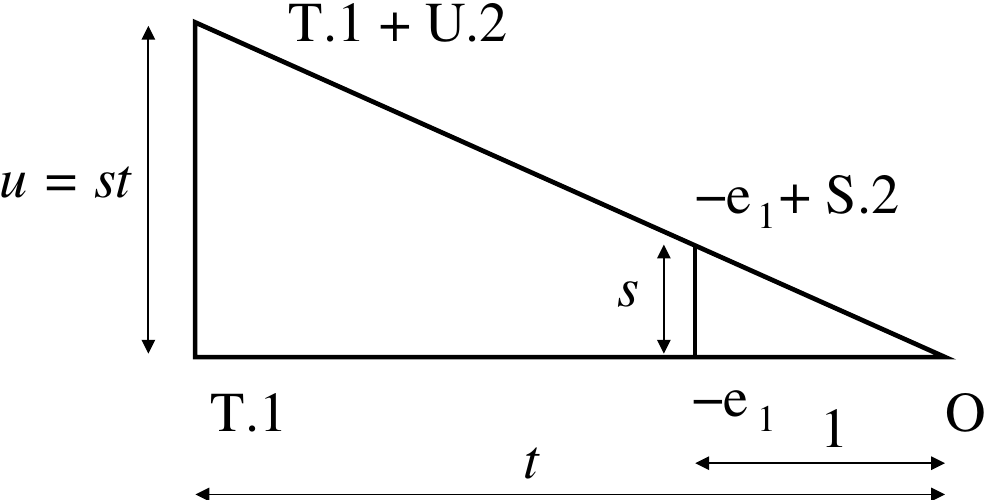}
\caption{Defining multiplication}
\label{fig:mult}
\end{center}
\end{figure}

Assuming $\Ve_1$ and $\Ve_2$ are points such that $x\Ve_1 + y\Ve_2$ is
rotund whenever $x \le 0 \le y$, I claim that under our representation
of real numbers as pairs of vectors, the
following formula represents the graph of the multiplication function
restricted to non-negative operands:
\[
\begin{array}{lcl}
\pNNMult{\VS, \VT, \VU} &\IsDef&
\pOK{\VS} \And \pOK{\VT} \And \pOK{\VU} \And \VS \ge \VZero \And \VT \ge \VZero \\
&\And&
\pSD{-\Ve_1+\VS.2, \VT.1 + \VU.2}
\end{array}
\]
\noindent
For $\pNNMult{\VS, \VT, \VU}$ holds iff {\em(i)} $\VS$, $\VT$ and $\VU$ represent
numbers $s$, $t$ and $u$ with $s, t \ge 0$ and 
{\em(ii)} the triangles
$\bigtriangleup\VZero(-\Ve_1)(-\Ve_1+\VS.2)$ and $\bigtriangleup\VZero(\VT.1)(\VT.1+\VU.2)$ have corresponding edges parallel and so are similar, whence:
\[
u = \Norm{\VU.2} = \Norm{\VS.2}\frac{\Norm{\VT.1}}{\Norm{-\Ve_1}}
= st.
\]
(See figure~\ref{fig:mult} and note that by assumption all the
points labelled in the figure are rotund so that the final conjunct
in $\pNNMult{\VS, \VT, \VU}$ ensures that $-\Ve_1+\VS.2$ lies on 
the ray through the origin in the direction of $\VT.1+\VU.2$).

We extend our definition of multiplication to
arbitrary real operands:
\[
\begin{array}{lcr@{}c@{}r@{}c@{}l@{}r@{}r@{}r@{}l}
\pMultX{\VS, \VT, \VU} &\IsDef& 
   \VS \ge \VZero &{}\And{}& \VT \ge \VZero &{}\And{}&
        \pNNMult{&\VS, &\VT, &\VU} \\
&\Or& \VS < \VZero &{}\And{}& \VT \ge \VZero &{}\And{}&
        \pNNMult{&-\VS, &\VT, &-\VU} \\
&\Or& \VS \ge \VZero &{}\And{}& \VT < \VZero &{}\And{}&
        \pNNMult{&\VS, &-\VT, &-\VU} \\
&\Or& \VS < \VZero &{}\And{}& \VT < \VZero &{}\And{}&
        \pNNMult{&-\VS, &-\VT, &\VU}&. \\
\end{array}
\]
\noindent The additive version of $\pG{s, t}$ may now be given as follows:
\[
\begin{array}{lcl}
\pG{\VS, \VT, \VU_1} &\IsDef&
\VS > \VZero \And \VT > \VZero \And
\pMultX{\VOne  + \VS, \VT, \VU_1} \\
&\And&
\Norm{(\VOne + \VT).1 + \VU_1.2} = 1 + \Norm{\VS.2} +\Norm{\VU_1.2}  \\
\end{array}
\]
\noindent so that $\pG{\VS, \VT, \VU_1}$ holds iff
$\VS$, $\VT$ and $\VU_1$ represent positive
real numbers $s$, $t$ and $u_1$ such that
$u_1 = (1+s)t$
and $\Norm{{-}(1+t)\Ve_1 + (1 + s)t\Ve_2} = 1 + s + u_1 = (1+s)(1+t)$,
in which case the design of the unit circle in $\sLI$
ensures that $t = 2s + s^2 + \frac{1}{M}\Sin{s}$
where $M$ is a positive integer constant. (See \SAHr{lemma~43}
and the discussion of the formulas $\pGamma{x, y}$ and
$\pG{s, t}$ that precedes it).

For the additive version of $\pSIN{s, t}$ we define:
\[
\begin{array}{lcl}
\pSIN{\VS, \VT, \VU_1, \VU_2} &\IsDef& 
\pG{\VS, 2\VS + \VU_2 + \frac{1}{M}\VT, \VU_1} \And \pMultX{\VS, \VS, \VU_2}
\end{array}
\]
\noindent so that $\pSIN{\VS, \VT, \VU_1, \VU_2}$ holds iff
$\VS$, $\VT$, $\VU_1$ and $\VU_2$ represent
real numbers $s$, $t$, $u_1$ and $u_2$ such that $s > 0$,
$u_1 = (1+s)(2s + s^2 + \frac{t}{M})$, $u_2 = s^2$ and $t = \Sin{s}$.

The periodicity property can be expressed as follows
with $\VA$, $\VS$ and $\VT$ representing the numbers $a$, $s$ and $t$
of the non-additive formulation and $\VV_1, \ldots \VV_5$ 
representing the results of various intermediate calculations:
\[
\begin{array}{lcl}
\pPeriodic &\IsDef& 
\pOK{\VA} \And \VA > \VZero
\\
&\And&
 (\VZero < \VS < 2\VA \And \pSIN{\VS, \VT, \VV_1, \VV_2} \Imp (\VT = \VZero \Iff \VS = \VA)) \\
&\And&
 (\pSIN{\VS, \VT, \VV_1, \VV_2} \And \pSIN{\VS + \VA, \VV_3, \VV_4, \VV_5} \Imp
	\VV_3 = -\VT)
\end{array}
\]
Assume that $\Ve_1$, $\Ve_2$, $\VA.1$ and $\VA.2$ have been interpreted in some normed
space and that $\pDef$ and $\pPeriodic$ hold under all extensions of this
interpretation to the other free variables of those formulas. The formula
$\ex{\VT\,\VU_1\,\VU_2} \pSIN{\VS, \VT, \VU_1, \VU_2} \And \VT = \VZero$ will
then hold iff $\VS$ is a positive integer multiple of $\VA$.

There follow the additive versions of the formula $\Func{N}$
to characterize the representations of natural numbers
and of the formula $\sf\Pi$ to characterize the representation
of the number $\pi$.
\[
\begin{array}{rcl}
\pN{\VX, \VU_1, \VU_2, \VU_3} &\IsDef& \pSIN{\VU_1, \VZero, \VU_2, \VU_3} \And \pMultX{\VX + \VOne, \VA, \VU_1}\\
\pPi{\VX, \VU_1, \VU_2} &\IsDef& \VX < \VFour  \And \pSIN{\VX, \VZero, \VU_1, \VU_2}.
\end{array}
\]
We now have all we need to give the additive version of the sentence $A$ in the proof
of \SAHr{theorem~45}. $A$ will hold in $\sLI$ and will characterize
a class of spaces in which for
suitable interpretations of $\Ve_1$, $\Ve_2$ and $\VA$,
$\ex{\VU_1\,\VU_2\,\VU_3}\pN{\VX, \VU_1, \VU_2, \VU_3}$
defines precisely those $\VX$ that represent natural numbers.
Writing $\VA$, $\VU_1$ \ldots between $\forall$ and $\cdot$ to stand for quantification over
$\VA.1$, $\VA.2$, $\VU_1.1$, $\VU_1.2$ \ldots, we may then define $A$ as follows:
\[
\begin{array}{lcl}
A &\IsDef&
\all{\Ve_1\,\Ve_2\,\Vw_1\,\Vw_2\,\Vw_3\,\VA\,\VU_1\,\VU_2\,\Vx\,\Vy\,\Vz\,\VS\,\VT\,\VV_1\,\VV_2\,\VV_3\,\VV_4\,\VV_5} \\
&    & \pW{\Ve_1, \Ve_2, \Vw_1, \Vw_2, \Vw_3} \And \pPi{\VA, \VU_1, \VU_2} \Imp
       \pDef \And \pPeriodic
\end{array}
\]
\noindent
The formula $\pDef$ here ensures that if $A$ and
$\pW{\Ve_1, \Ve_2, \Vw_1, \Vw_2, \Vw_3}$ hold
and if $\pPi{\VA, \VU_1, \VU_2}$ is satisfiable, then
the point $x\Ve_1 + y\Ve_2$ is rotund if $x \le 0 \le y$.
Thus $\pMultX{\VS, \VT, \VU}$ represents
multiplication under any interpretation
of $\Ve_1$, $\Ve_2$ that can be extended to make $\pW{\Ve_1, \Ve_2, \Vw_1, \Vw_2, \Vw_3}$ and  $\pPi{\VA, \VU_1, \VU_2}$ hold.
$\pPeriodic$ then implies that
$\ex{\VU_1\,\VU_2\,\VU_3}\pN{\VX, \VU_1, \VU_2, \VU_3}$
holds iff $\VX = (-x\Ve_1, x\Ve_2)$ for some $x \in \nat$.

Now let $Q(x_1, \ldots, x_k)$ be a quantifier-free formula in the language of
arithmetic. We can transform $Q(x_1, \ldots, x_k)$ into a logically equivalent
formula with the same free variables and having the following form \[
\ex{s_1 \ldots s_m\,t_1 \ldots t_m\,z_1 \ldots z_m} z_1 = s_1t_1 \And \ldots \And z_m = s_mt_m \And Q_1
\]
\noindent 
where the formula $Q_1$ is additive and quantifier-free.
Define $B$ as follows:
\[
\begin{array}{lcc}
B &\IsDef&
\all{\Ve_1\,\Ve_2\,\Vw_1\,\Vw_2\,\Vw_3\,\VA\,\VU_1\,\VU_2\,
   \VS_1\ldots\VS_{4m}\,\VT_1\ldots\VT_{4m}\,
   \VZ_1\ldots\VZ_{m} \\
&    &     s_1 \ldots s_m\,t_1 \ldots t_m\,z_1 \ldots z_m\,
                \VX_1\ldots\VX_{4k}} \hspace{1.1in}\\
&    & \pW{\Ve_1, \Ve_2, \Vw_1, \Vw_2, \Vw_3} \And \pPi{\VA, \VU_1, \VU_2} \\
&\And&     {\displaystyle \Ands_{i=1}^{m}}
	\left(
        \begin{array}{@{}c@{}c}
                &
                    \pN{\VS_i, \VS_{m+i}, \VS_{2m+i}, \VS_{3m+i}} \\
            {} \And {} &
                           \pN{\VT_i, \VT_{m+i}, \VT_{2m+i}, \VT_{3m+i}}  \\
            {} \And {} &
                          \pMultX{\VS_i, \VT_i, \VZ_i} \\
            {} \And {} &
                             s_i = \Norm{\VS_i.2} \And
                             t_i = \Norm{\VT_i.2} \And
                             z_i = \Norm{\VZ_i.2}
        \end{array}\right) \\
&\And&     {\displaystyle\Ands_{i=1}^{k}}
                (\pN{\VX_i, \VX_{k+i}, \VX_{2k+i}, \VX_{3k+i}} \And
                  x_i = \Norm{\VX_i.2}) \\
&\Imp& \Not Q_1
\end{array}
\]
\noindent In $\sLI$, $A$ holds and the formulas $\pW{\Ve_1, \Ve_2, \Vw_1, \Vw_2,
\Vw_3}$ and $\pPi{\VA, \VU_1, \VU_2}$ are simultaneously satisfiable: in such a
space, $B$ and hence $A \Imp B$ hold iff $Q_1$ and hence also $Q(x_1, \ldots, x_k)$ are unsatisfiable
in $\nat$.  $A \Imp B$ clearly holds in a space in which $A$ does not hold or
in which $\pW{\Ve_1, \Ve_2, \Vw_1, \Vw_2, \Vw_3}$ and  $\pPi{\VA, \VU_1,
\VU_2}$ are not simultaneously satisfiable.  So $A \Imp B$ is valid in a class
of spaces including $\sLI$ iff $Q(x_1, \ldots, x_k)$ is unsatisfiable in
$\nat$.  Thus a decision procedure for valid $\fAIA$ sentences could decide $A
\Imp B$ and so give a decision procedure for unsatisfiable quantifier-free
formulas of arithmetic, which is impossible.  This gives the case $d = 2$ of
the following theorem:

\begin{Theorem}\label{thm:ns-additive-aia-valid-undec}
Let $d \in \{2, 3, 4, \ldots\} \cup \{\infty\}$ and
let $\cal M$ be any class of normed spaces that includes all Banach
spaces of dimension $d$.
The set of additive $\fAIA$ sentences that are valid in $\cal M$ is
undecidable.
\end{Theorem}
\Proof
We have done the case $d = 2$.
For $d > 2$,  let $V$ be a Hilbert space of dimension $d - 2$, let
$W$ be the 2-sum of $\sLI$ and $V$, and identify $\sLI$ with the subspace
$\sLI \times 0$ of $W$. See Appendix~\ref{sec:two-sum} for the definition
of the 2-sum and a proof that if a line segment $[\Vu, \Vv]$
lies in the unit sphere of $W$, then it is parallel to $\sLI$.

Assume $\pW{\Ve_1, \Ve_2, \Vw_1, \Vw_2, \Vw_3}$ holds in $W$. Then
$[\Vw_2, \Vw_3]$ is a straight line segment in the unit sphere of $W$
of the maximal length $2r$ and hence is contained in $\sLI$.
So $[\Vw_1, \Vw_2]$ is a line segment in the unit sphere with one endpoint
in $\sLI$ and hence is also contained in $\sLI$.
However, as things stand $\Ve_1$ and $\Ve_2$ need not lie in $\sLI$.
To avoid this problem, modify $A$ and $B$ to
give new sentences $A'$ and $B'$ by
conjoining the following formula (*) to $\pW{\Ve_1, \Ve_2, \Vw_1, \Vw_2, \Vw_3}$
and universally quantifying over $\Va_1$, $\Va_2$, $\Vb_1$ and $\Vb_2$.
\[
\begin{array}{ccc}
\hspace{.75in}~&
\begin{array}{cc}
    & \Ve_1 = \Va_1 + \Vb_1 \And \pPar{\Va_1, \Vw_1} \And \pPar{\Vb_1, \Vw_2} \\
\And& \Ve_2 = \Va_2 + \Vb_2 \And \pPar{\Va_2, \Vw_1} \And \pPar{\Vb_2, \Vw_2} 
\end{array}
    \hspace{.75in}\mbox{(*)}
\end{array}
\]
\noindent
Even though the $\Vw_j$ are not rotund,
if $\Vv \in W$ satisfies $\pPar{\Vv, \Vw_j}$,
then one of the line segments $[\Vv, \pm\frac{\Norm{\Vv}}{\Norm{\Vw_j}}\Vw_j]$
is contained in the $\sLI$-circle $S_{\Norm{\Vv}}$  and so $\Vv \in \sLI$.
Hence (*) implies that $\{\Ve_1, \Ve_2\} \subseteq \sLI$.
Conversely, (*) holds in $\sLI$ if we put
$\Ve_i = x_i\Vw_1 + y_i\Vw_2$ and $\Va_i = x_i\Vw_1$ and
$\Vb_i = y_i\Vw_2$.
The argument for $A$, $B$ and $\sLI$ may now be applied to $A'$, $B'$ and $W$.
\Done

Via very simple reductions, $\fAIA$ formula have both $\fAE$ and $\fEA$
equivalents.
Hence the theorem has the following immediate corollaries:
\begin{Corollary}\label{thm:ns-additive-ae-valid-undec}
Let $d \in \{2, 3, 4, \ldots\} \cup \{\infty\}$ and
let $\cal M$ be any class of normed spaces that includes all Banach
spaces of dimension $d$.
The set of additive $\fAE$ sentences that are valid in $\cal M$ is
undecidable. \Done
\end{Corollary}

\begin{Corollary}\label{thm:ns-additive-ea-valid-undec}
Let $d \in \{2, 3, 4, \ldots\} \cup \{\infty\}$ and
let $\cal M$ be any class of normed spaces that includes all Banach
spaces of dimension $d$.
The set of additive $\fEA$ sentences that are valid in $\cal M$ is
undecidable. \Done
\end{Corollary}

\section*{Acknowledgments}
I am indebted to Robert M. Solovay for finding a number of infelicities in an
earlier draft of this paper. Any errors that remain are, of course, mine.
\bibliographystyle{plain}

\bibliography{bookspapers}

\newpage
\appendix
\section{Intersecting Circles}~\label{sec:intersecting-circles}

We will prove lemma~\ref{lma:two-point-lemma} using a general fact about
the intersection of the boundaries of two convex bodies in the plane.
If $X$ and $Y$ are convex bodies in $\real^{n}$ then the intersection
$\Bndry{X} \cap \Bndry{Y}$ of their boundaries can be topologically very
complicated: e.g., if $C$ is any closed subset of the boundary $S^{n-1}
= \Bndry{D^{n}}$ of the unit ball $D^{n}\subseteq \real^{n}$ and $X$ is
the convex hull of $C$, then $\Bndry{X} \cap S^{n-1} = C$.  The
situation when $Y$ is obtained from $X$ by scaling and translation is
much simpler. Here, we are only concerned with the case $n = 2$.  For
background on the following theorem see 
\cite[Section 3.3]{martini-et-al01}. As the literature on
this is not very accessible, we give a proof here based on a
construction of Sch\"{a}ffer \cite[Lemma 4.3]{schaeffer67}.
\begin{Theorem}\label{thm:intersecting-circles}
Let $D_1 \subseteq \real^2$ be a convex body, let
$c \in \realgtz$ and let $\Vv \in \real^2$.
Define $f : \real^2 \To \real^2$ by $f(\Vx) = c\Vx + \Vv$,
let $D_2 = f(D_1)$ and let $S_i = \Bndry{D_i}$, $i = 1, 2$.
Then one of the following holds:
\begin{enumerate}
\item \label{case:disjoint}
$S_1 \cap S_2 = \emptyset$;
\item \label{case:eq}
$S_1 \cap S_2 = S_1 = S_2$;
\item \label{case:one-comp}
$S_1 \cap S_2$ has one connected component which is either
{\em(a)} a point or a proper closed line segment, or
{\em(b)} the union of two proper closed
line segments that meet at a common endpoint and are not collinear.
\item \label{case:two-comps}
$S_1 \cap S_2$ has two connected components each of which is
either a point or a proper closed line segment.
\end{enumerate}
\end{Theorem}
\Proof
Case {\em(i)} $\Int{D_1} \cap \Int{D_2} = \emptyset$:
(here, by definition, $\Int{D_{i}} = D_i \Diff \Bndry{D_i}$.) If
$\Vx \in S_{2-i} \cap \Int{D_i}$, then for $\Vy \in \Int{D_{2-i}}$
sufficiently close to $\Vx$, we have
$\Vy \in \Int{D_1} \cap \Int{D_2}$, contradicting our assumption.
So, $S_1 \cap S_2 = D_1 \cap D_2$,
which is a closed bounded convex subset of $\real^2$ with empty interior.
Such a subset is either empty giving outcome~\ref{case:disjoint}
or a point or a proper closed line segment giving outcome~\ref{case:one-comp}{\em(a)}.

Case {\em(ii)} $\Int{D_1} \cap \Int{D_2} \not= \emptyset \And c = 1$:
in this case, $f$ is translation by $\Vv$.
If $\Vv = \VZero$, then $f$ is the identity function, so
that $S_1 = S_2$ giving us outcome~\ref{case:eq}.
So we may assume that $\Vv \not= \VZero$.
Let $\Vw \not= \VZero$ be orthogonal to $\Vv$ and for any $y \in \real$,
let $l(y)$ be the line through $y\Vw$ parallel to $\Vv$.
Let $X = \{y \in \real \ST l(y) \cap \Int{D_1} \cap \Int{D_2} \not= \emptyset\}$.
Then $X$, is easily seen to be a non-empty open, connected and bounded
subset of $\real$,,
i.e., $X$ is the open interval $(y_1, y_2)$ for some $y_1 < y_2$.
It is also easy to see that $D_1 \cap D_2$ is contained in the closed strip
bounded by $l(y_1)$ and $l(y_2)$ and meets both $l(y_1)$ and $l(y_2)$.  If $y_1
< y < y_2$, then $l(y) \cap \Int{D_1} \not= \emptyset \not= l(y) \cap
\Int{D_2}$, so that there are vectors $\Va$ and $\Vb$ such that $l(y) \cap D_1
= [\Va, \Vb]$, $l(y) \cap D_2 = [\Va + \Vv, \Vb + \Vv]$ and $(\Va, \Vb) \cap
(\Va + \Vv, \Vb + \Vv) \not= \emptyset$. Hence $l(y) \cap S_1 = \{\Va, \Vb\}$
cannot meet $l(y) \cap S_2 = \{\Va + \Vv, \Vb + \Vv\}$.  It follows that $S_1
\cap S_2 = (l(y_1) \cap D_1 \cap D_2) \cup (l(y_2) \cap D_1 \cap D_2)$.  Each
of the two sets $l(y_i) \cap D_1 \cap D_2$ ($i \in \{1, 2\}$) is a closed
bounded convex subset of $l(y_i)$, i.e., is empty, a point or a line segment,
giving outcome~\ref{case:one-comp}{\em(a)} or outcome~\ref{case:two-comps} according as
one or both of these sets is non-empty.

Case {\em(iii)} $\Int{D_1} \cap \Int{D_2} \not= \emptyset \And c \not= 1$:
in this case, $f$ has a unique fixed point $\Vz = \Vv/(1-c)$.
As the property we wish to prove is invariant under translations,
we may arrange for $\Vz = \Vv = \VZero$, by translating both
$D_1$ and $D_2$ by $-\Vz$.
So we may assume that $f(\Vx) = c\Vx$, i.e., $f$ is scaling by $c$.
If $\Vx \in (\Int{D_1} \cap \Int{D_2}) \Diff \{\VZero\}$, there are $x, y \in
\real$ such that the line $\overline{\VZero\Vx}$ meets $D_1$
in the line segment $[x\Vx, y\Vx]$ and meets $D_2$ in the
line segment $[cx\Vx, cy\Vx]$, so that $l \cap S_1 = \{x\Vx, y\Vx\}$
and $l \cap S_2 = \{cx\Vx, cy\Vx\}$; moreover if $\VZero \in \Int{D_1}$
we must have $xy < 0$ and so $l \cap S_1 \cap S_2 = \emptyset$.
It follows that if $\VZero \in \Int{D_1}$, then we have
outcome~\ref{case:disjoint}.
If $\VZero \not\in \Int{D_1}$, pick $\Vx \in S_1$
with $\Norm{\Vx}$ minimal (possibly $\Vx = \VZero$) and let
$l$ be a line through $\VZero$ parallel to a supporting line for $D_1$
at $\Vx$; then $l \cap D_1 \subseteq \{\VZero\}$ and hence
$l \cap D_2 = f(l \cap D_1) =l \cap D_1$.
Thus $D_1 \cup D_2$ lies in one of the closed halfspaces determined by
$l$ and we may choose $\Vw \in l \Diff \{\VZero\}$ such that
$0 \le \angle \Vx\Vw \le \pi$ for every $\Vx \in (D_1 \cup D_2) \Diff \{\VZero\}$.
For $0 \le \theta \le \pi$, let $r(\theta)$ be the ray through $\VZero$
such that $\angle r(\theta)\Vw = \theta$ and let
$Y = \{\theta \in [0, \pi] \ST r(\theta) \cap \Int{D_1} \cap \Int{D_2} \not=
\emptyset\}$. Then $Y$
is non-empty by assumption and is easily seen to be open, connected and bounded,
i.e., $Y$ is the open interval $(\theta_1, \theta_2)$ for some $\theta_1 < \theta_2$.
It is easy to verify that $D_1 \cap D_2$ is contained in the closed cone bounded by $r(\theta_1)$
and $r(\theta_2)$ and meets both $r(\theta_1)$ and $r(\theta_2)$.
Now if $\theta_1 < \theta < \theta_2$ and $\Vx \in r(\theta) \Diff
\{\VZero\}$, then there are $x, y$ with $0 \le x < y$, such that
$r(\theta) \cap D_1 = [x\Vx, y\Vx]$ and $r(\theta) \cap D_2 = [cx\Vx,
cy\Vx]$. So, as $\theta_1 < \theta < \theta_2$,
we must have $(x, y) \cap (cx, cy) \not= \emptyset$, implying
that $r(\theta) \cap S_1 \cap S_2 \subseteq \{\VZero\} \subseteq
r(\theta_1) \cup r(\theta_2)$. 
Since $D_1 \cup D_2$ is contained in the union of the $r(\theta)$
with $\theta \in [\theta_1, \theta_2]$, it follows that
$S_1 \cap S_2 = (r(\theta_1) \cap D_1 \cap D_2) \cup (r(\theta_2) \cap D_1
\cap D_2)$.
Each of the two sets $r(\theta_i) \cap D_1 \cap D_2$
is a convex subset of $r(\theta_i)$, i.e., is empty, a point or a proper line
segment. If just one of these sets is empty, we have
outcome~\ref{case:one-comp}{\em(a)}, while, if both are non-empty, we have
outcome~\ref{case:two-comps} if they are disjoint and
outcome~\ref{case:one-comp}{\em(a)} or~\ref{case:one-comp}{\em(b)} if
they meet.  
\Done

We can now prove lemma~\ref{lma:two-point-lemma}:
writing $D_V(\Vp, r)$ for the closed disc of radius $r$ about $\Vp$ in $V$,
apply theorem~\ref{thm:intersecting-circles} with
$D_1 = D_V(\Vp, r)$, $c = \frac{s}{r}$ and $\Vv = \Vq -\frac{s}{r}\Vp$,
so that $D_2 = D_V(\Vq, s)$,
$S_1 = \Bndry{D_V(\Vp, r)} = S_V(\Vp, r)$ and
$S_2 = \Bndry{D_V(\Vq, s)} = S_V(\Vq, s)$.
By assumption $\Va \in S_1 \cap S_2$ and at least one
of $\Va - \Vp$ and $\Va - \Vq$ is rotund so that $\Va$
cannot belong to any proper line segment in $S_1 \cap S_2$.
Similarly, $\Vb \in S_1 \cap S_2$ and $\Vb$ cannot belong to any proper line
segment in $S_1 \cap S_2$.  As $\Vp \not= \Vq$, $S_1 \not= S_2$, and so, as
$\Va \not= \Vb$, the only possibility in the conclusion of
theorem~\ref{thm:intersecting-circles} is that $S_1 \cap S_2$ has two connected
components and these are points. So $S_1 \cap S_2 = \{\Va,
\Vb\}$ as the lemma claims.

\section{The unit disc in a 2-sum}\label{sec:two-sum}
If $U$ and $V$ are normed spaces
their 2-sum $U \mathop{{\times}_2} V$ is the product vector space $U \times V$ equipped
with the norm defined by:
$$
\Norm{(\Vu, \Vv)} = \sqrt{\Norm{\Vu}_U^2 + \Norm{\Vv}_V)^2}
$$
Writing $\Norm{\_}_e$ for the euclidean norm on $\real^2$, we
may also express this as
$\Norm{(\Vu, \Vv)} = \Norm{(\Norm{\Vu}_U, \Norm{\Vv}_V)}_e$,
a point of view which is helpful in proving the following
theorem about the 2-sum of a normed space and an inner product space. This
has been used in this note in the proof of
theorem~\ref{thm:ns-additive-aia-valid-undec} and in the proof of
\SAHr{theorem~41}.

\begin{Theorem}\label{thm:normed-times-hilbert}
Let $U$ be a normed space, let $V$ be an inner product space
and let $W = U \mathop{{\times}_2} V$ be their 2-sum.
If $[\Vw_1, \Vw_2]$ is a line segment contained in the unit
sphere of $W$, then $[\Vw_1, \Vw_2]$ is parallel to $U$,
\end{Theorem}
\Proof
Writing $\Vw_1 = (\Vu_1, \Vv_1)$ and $\Vw_2 = (\Vu_2, \Vv_2)$,
we have to prove that $\Vv_1 = \Vv_2$.
Let us write  $\Vw_3 = (\Vu_3, \Vv_3) = \Frac{\Vw_1 + \Vw_2}{2}$ so that
$\Vu_3  = \Frac{\Vu_1 + \Vu_2}{2}$ and $\Vv_3 = \Frac{\Vv_1 + \Vv_2}{2}$.
For $i = 1, 2, 3$, let us write $x_i = \Norm{\Vu_i}_U$, $y_i = \Norm{\Vv_i}_V$
and $z_i = \Norm{\Vw_i}_W$
so that $z_i = \Norm{(x_i, y_i)}_e = 1$.
I claim that $y_1 = y_2 = y_3$ and
this will complete the proof since it implies that
$[\Vv_1, \Vv_2] \subseteq S_{\Vv_1} = 
   \{\Vv : V \ST \Norm{\Vv}_V = \Norm{\Vv_1}_V\}$,
but then, as $V$ is an inner product space, $S_{\Vv_1}$ contains no proper
line segments and we can only have $\Vv_1 = \Vv_2$.
To prove the claim, first note that,
using the triangle inequality in $U$ and $V$, we have:
$$
\begin{array}{rclcl}
0 &\le& x_3 &\le& \FRAC{x_1 + x_2}{2} \\
0 &\le& y_3 &\le& \FRAC{y_1 + y_2}{2}
\end{array}
$$
\noindent
Using the above inequalities and the triangle inequality for $\Norm{\_}_e$
we find:
$$
\begin{array}{rcl}
1 &=& \sqrt{x_3^2 + y_3^2} \\
  &\le&
    \sqrt{(\FRAC{x_1 + y_1}{2})^2 + (\FRAC{x_2 + y_2}{2})^2} =
      \Norm{\FRAC{(x_1, x_2) + (y_1, y_2)}{2}}_e\\
  &\le&
    \Norm{\FRAC{(x_1, x_2)}{2}}_e + \Norm{\FRAC{(y_1, y_2)}{2}}_e
        = \FRAC{z_1}{2}+ \FRAC{z_2}{2}= 1
\end{array}
$$
\noindent so $\Norm{\Frac{(x_1, x_2) + (y_1, y_2)}{2}}_e = 1$
whence
the line segment $[(x_1, y_1), (x_2, y_2)]$ is contained
in the euclidean unit circle in $\real^2$ and so
must reduce to a point, implying
that $y_1 = y_2 = y_3$ and completing the proof of the claim and
hence the theorem.
\Done
\typeout{==== REMEMBER TO CHECK THE DATE ===}
\end{document}